\documentclass[12pt]{article}
\usepackage{amsmath,amssymb,latexsym,epsfig}
\newcommand{\beq}[1]{\begin{equation}\label{#1}}
\newcommand{\eeq}{\end{equation}}

\newcommand{\set}[1]{\left\{#1\right\}}

\newcommand{\ignore}[1]{}
\def\qs{{\bf qs}}

\def\cE{{\cal E}}

\def\cR{{\cal R}}

\def\E{{\bf E}}

\def\cP{{\cal P}}

\def\a{\alpha}

\def\d{\delta}
\def\D{\Delta}
\def\e{\epsilon}
\def\f{\phi}

\def\G{\Gamma}

\def\th{\theta}

\def\n{\nu}
\def\p{\pi}

\def\r{\rho}

\def\s{\sigma}

\def\t{\tau}

\def\Pr{\mbox{{\bf Pr}}}

\def\whp{{\bf whp}}

\newtheorem{lemma}{Lemma}
\newtheorem{theorem}{Theorem}

\newcommand{\brac}[1]{\left( #1\right)}
\newcommand{\bfrac}[2]{\brac{\frac{#1}{#2}}}

\newcommand{\proofstart}{{\bf Proof\hspace{2em}}}

\newcommand{\proofend}{\hspace*{\fill}\mbox{$\Box$}}

\newcommand{\rdup}[1]{\left\lceil #1 \right\rceil }
\newcommand{\rdown}[1]{\left\lfloor #1 \right\rfloor }

\newtheorem{remark}{Remark}
\begin{document}
\title{{\bf Packing Hamilton Cycles in Random and Pseudo-Random Hypergraphs}}
\author{Alan Frieze\thanks{
Department of Mathematical Sciences, Carnegie Mellon University,
Pittsburgh PA15213, U.S.A. Supported in part by NSF grant
DMS-0753472.} \and Michael Krivelevich\thanks{ School of
Mathematical Sciences, Raymond and Beverly Sackler Faculty of Exact
Sciences, Tel Aviv University, Tel Aviv, 69978, Israel. Email:
krivelev@post.tau.ac.il. Research supported in part by a USA-Israel
BSF grant, by a grant from the Israel Science Foundation and by a
Pazy Memorial Award.} }
\maketitle
\begin{abstract}
We say that a $k$-uniform hypergraph $C$ is a Hamilton cycle of type
$\ell$, for some $1\le \ell \le k$, if there exists a cyclic
ordering of the vertices of $C$ such that every edge consists of $k$
consecutive vertices and for every pair of consecutive edges
$E_{i-1},E_i$ in $C$ (in the natural ordering of the edges) we have
$|E_{i-1}-E_i|=\ell$. We prove that for $\ell \le k\le 2\ell$, with
high probability almost all edges of a random $k$-uniform hypergraph
$H(n,p,k)$ with $p(n)\gg \log^2 n/n$ can be decomposed into edge
disjoint type $\ell$ Hamilton cycles. We also provide sufficient
conditions for decomposing almost all edges of a pseudo-random
$k$-uniform hypergraph into type $\ell$ Hamilton cycles, for $\ell
\le k\le 2\ell$. For the case $\ell=k$ these results show that
almost all edges of corresponding random and pseudo-random
hypergraphs can be packed into disjoint perfect matchings.
\end{abstract}

\section{Introduction}
The subject of Hamilton graphs and Hamiltonicity-related problems is
undoubtedly one of the most central in Graph Theory, with great many
deep and beautiful results obtained. Hamiltonicity problems occupy a
place of honor in the theory of random graphs too, the reader can
consult the monographs of Bollob\'as \cite{Bol} and of Janson, \L
uczak and Ruci\'nski \cite{JLR} for an account of some of the most
important results related to Hamilton cycles in random graphs. Of 
particular relevance to the current work is a previous result of the
authors \cite{FK1} who proved that for edge probability $p=p(n)\ge
n^{-\e}$ for some constant $\e>0$, \whp\footnote{An event $\cE_n$
occurs {\em with high probability}, or \whp\ for brevity, if
$\lim_{n\rightarrow\infty}\Pr(\cE_n)=1$.} almost all edges of the
random graph $G(n,p)$ can be packed into edge disjoint Hamilton
cycles.

Quite a few results about Hamiltonicity of pseudo-random graphs are
available too. Informally, a graph $G=(V,\cE)$ with $|V|=n$ vertices
and $|\cE|=m$ edges is pseudo-random if its edge distribution is
similar, in some well defined quantitative way, to that of a truly
random graph $G(n,p)$ with the same expected density
$p=m\binom{n}{2}^{-1}$. A thorough discussion about pseudo-random
graphs, their alternative definitions and properties can be found in
survey \cite{KS}. It is well known that pseudo-randomness of graphs
can be guaranteed by imposing conditions on vertex degrees and
co-degrees (see, e.g., \cite{Thom}, \cite{CGW}); we will adopt a
similar approach later in the paper when discussing pseudo-random
hypergraphs. There are known sufficient criteria for Hamiltonicity
in pseudo-random graphs. Also, the above mentioned result of
\cite{FK1} can be extended to the pseudo-random case as well. Since
we will employ this result in our arguments, let us state it here
formally. A graph $G$ on vertex set $[n]$ is {\em $(\a,\e)$-regular}
if
\begin{description}
\item[$Q_a$:] $\d(G)\geq (\a-\e) n$.
\item[$Q_b$:] If $S,T$ are disjoint subsets of $[n]$ and $|S|,|T|\geq \e
n$ then $\left|\frac{e_G(S,T)}{|S|\,|T|}-\a\right|\leq \e$, where
$e_G(S,T)$ is the number of $S-T$ edges in $G$.
\end{description}
The following is implied by the main theorem of \cite{FK1}:
\begin{theorem}\label{FK1}
Let $G$ be an $(\a,\e)$-regular graph with $n$ vertices where
$$\a\gg\e\ and\ \a\e^3\gg\frac{1}{(n\log n)^{1/2}}.$$
Then $G$ contains at least $(\a/2-4\e)n$ edge disjoint Hamilton
cycles.
\end{theorem}
\begin{remark}\label{rem1}
Theorem 2 of \cite{FK1} only claims to be true for $\a$ constant.
This was an unfortunate over-cautious statement. The real condition
should be the one given in the above theorem.
\end{remark}

In contrast, much less is known about Hamiltonicity in hypergraphs
in general and in random and pseudo-random hypergraphs in
particular. Formally, a {\em hypergraph} $H$ is an ordered pair
$H=(V,\cE)$, where $V$ is a set of vertices, and $\cE$ is a family
of distinct subsets of $V$, called edges. A hypergraph $H$ is {\em
$k$-uniform} if all edges of $H$ are of size $k$. It is generally
believed that $k$-uniform hypergraphs for $k\ge 3$  are much more
complicated objects of study than graphs (corresponding to $k=2$).
Specifically for Hamiltonicity, even extending the definition of a
Hamilton cycle in graphs to the case of (uniform) hypergraphs is not
a straightforward task. In fact, several alternative definitions are
possible. In this paper (in some departure from a relatively
standard notation) we will use the following definition. Denote
$$\n_i=\frac{n}{i},\,1\leq i\leq k.$$
Suppose that $1\leq \ell\leq k$. A {\em type $\ell$ Hamilton cycle}
in a $k$-uniform hypergraph $H=(V,\cE)$ on $n$ vertices is a
collection of $\n_\ell$ edges of $H$ such that for some cyclic order
of $[n]$ every edge consists of $k$ consecutive vertices and for
every pair of consecutive edges $E_{i-1},E_i$ in $C$ (in the natural
ordering of the edges) we have $|E_{i-1}\setminus E_i|=\ell$. Thus, in a type
$\ell$ Hamilton cycle the sets $C_i=E_i\setminus
E_{i-1},\,i=1,2,\ldots,\n_\ell$, are a partition of $V$ into sets of
size $\ell$. (An obvious necessary condition for the existence of a
cycle of type $\ell$ in a hypergraph on $n$ vertices is that $\ell$
divides $n$. We thus always assume, when discussing Hamilton cycles
of type $\ell$, that this necessary condition is fulfilled.) In the
literature, when $\ell=1$ we have a {\em tight} Hamilton cycle and
when $\ell=k-1$ we have a {\em loose} Hamilton cycle. In the extreme
case $\ell=k$ the notion reduces to that of a perfect matching in a
hypergraph.

Several recent papers (see, e.g., \cite{HS}, \cite{KKMO},
\cite{KMO}) provided sufficient conditions for the existence of a
type $\ell$ Hamilton cycle in a $k$-uniform hypergraph $H$ on $n$
vertices in terms of the minimum number of edges of $H$ passing
through any subset of $k-1$ vertices, thus extending the classical
Dirac sufficient condition for graph Hamiltonicity to the
hypergraph case. These results however appear to be of rather
limited relevance to the current paper, as here we are mostly
concerned with sparse hypergraphs (with $o(|V|^k)$ edges), while the
above mentioned results are for the (very) dense case.

The main goal of this paper at large is to study Hamiltonicity in
random and pseudo-random hypergraphs. A {\em random $k$-uniform
hypergraph} $H(n,p,k)$ is a hypergraph with vertex set
$\{1,\ldots,n\}=[n]$, where each $k$-tuple of $[n]$ is an edge of
the hypergraph independently with probability $p=p(n)$. For the case
$k=2$ the model $H(n,p,k)$ reduces to the classical binomial random
graph $G(n,p)$. Essentially nothing appears to be known about
Hamilton cycles in random hypergraphs. Even the most basic question
of the threshold for the appearance of a cycle of type $\ell$ in
$H(n,p,k)$ has not yet been addressed. One notable exception is the
case $\ell=k$, i.e., the case of perfect matchings -- a recent
striking result of Johannson, Kahn and Vu \cite{JKV} has established
the order of magnitude of the threshold for the appearance of a
perfect matching in a $k$-uniform random hypergraph.

In this paper, rather than studying the conditions for the existence of a single
Hamilton cycle, we study the conditions for the existence of a {\em packing of almost all edges} of a random
or a pseudo-random hypergraph into Hamilton cycles.
For $\ell \ge
k/2$ we manage to obtain non-trivial results in this direction.
 It appears that
the cases of small $\ell$ (where adjacent edges along the Hamilton
cycle have larger intersection) are harder. 

Our first result is about packing Hamilton cycles in random
hypergraphs.
\begin{theorem}\label{th5}
Suppose that $\ell\leq k\leq2\ell$ and suppose that
$np/\log^2n\to\infty$. Then \whp\ $H=H(n,p,k)$ contains a collection
of $(1-\e)\binom{n}{k}p/\n_\ell$ edge disjoint type $\ell$ Hamilton
cycles, where $\e=O((\log n/(np)^{1/2})^{1/2})=o(1)$.
\end{theorem}

Note that for the case $\ell=k$ the above theorem provides a
sufficient condition on the edge probability $p(n)$ for being able
to pack \whp\ almost all edges of $H(n,p,k)$ into perfect matchings.

Other results of the paper are about packing Hamilton cycles in
pseudo-random hypergraphs. For most part, we state the condition of
pseudo-randomness of a hypergraph in terms of the number of edges
through subsets of vertices of fixed size. These conditions are
suggested by the expected numbers of such edges in truly random
hypergraphs of the same edge density and are easily seen to hold
\whp\ in random hypergraphs. Thus our results about pseudo-random
hypergraphs are applicable to truly random instances as well.
Naturally, the direct approach of Theorem \ref{th5} provides a
better lower bound on the edge probability $p(n)$.

In this paper we are only able to deal with the case where $\ell\geq
k/2$. Let $H=([n],\cE)$ be a $k$-uniform
hypergraph with vertex set $[n]$ and $m$ edges. Its density $p=m/\binom{n}{k}$. For a set $X\subseteq [n]$ with $|X|=a<k$ we
define its neighbourhood $N_H(X)=\set{Y\in \binom{[n]}{k-a}: X\cup
Y\in \cE}$ and its degree $d_H(X)=|N_H(X)|$.

We first consider $k/2<\ell<k$ and list the following properties.
The value $\e$ will be a parameter of {\em regularity}.
\begin{description}
\item[$P_a$:] $\displaystyle{\min_{S\in \binom{[n]}{2(k-\ell)}}}\,d_H(S)\geq\rdup{(1-\e)\binom{n-2(k-\ell)}{2\ell-k}p}$.
\item[$P_b$:] $\displaystyle{\min_{S\in \binom{[n]}{2\ell-k}}}\,d_H(S)\geq \rdup{(1-\e)\binom{n-2\ell+k}{2(k-\ell)}p}$.
\item[$P_c$:] $\displaystyle{\max_{S\in \binom{[n]}{2(k-\ell)+1}}}\,d_H(S)\leq
 \rdown{(1+\e)\binom{n-2k+2\ell-1}{2\ell-k-1}p}$.
\item[$P_d$:] $\displaystyle{\max_{S\in \binom{[n]}{2\ell-k+1}}}\,d_H(S)\leq
 \rdown{(1+\e)\binom{n-2\ell+k-1}{2(k-\ell)-1}p}$.
\item[$P_e$:] $\displaystyle{\max_{\substack{S_1,S_2\in \binom{[n]}{2\ell-k}\\ S_1\cap S_2=\emptyset}}}\,|N_H(S_1)\cap
 N_H(S_2)|\leq \rdown{(1+\e)\binom{n-2(2\ell-k)}{2(k-\ell)}p^2}$.
\item[$P_f$:] $\displaystyle{\max_{\substack{S_1,S_2\in \binom{[n]}{2(k-\ell)}\\ |S_1\cap S_2|=0\ 
or\ k-\ell}}}\,|N_H(S_1)\cap
 N_H(S_2)|\leq \rdown{(1+\e)\binom{n}{2\ell-k)}p^2}$.
\end{description}
\begin{theorem}\label{th1}
Let $H=([n],\cE)$ be a $k$-uniform hypergraph with with $m$ edges,
and let $k/2<\ell< k$ and $1> \e^5\gg \log^3n/(n^{1/2}p^2).$ Suppose
that $H$ satisfies properties $\cP=\set{P_a,P_b,P_c,P_d,P_e,P_f}$. Then
$H$ contains a collection of $(1-2\e^{1/3})m/\n_\ell$ edge disjoint
type $\ell$ Hamilton cycles.
\end{theorem}
The restriction $1>\e$ is for relevance and the restriction $\e^5\gg
\log^3n/(n^{1/2}p^2)$ is used in the proof (see Lemma \ref{lem3}).\footnote{We use the notation
$a_n\gg b_n$ as shorthand for $a_n/b_n\to\infty$ as $n\to\infty$.}
The latter condition can be relaxed a little through a more careful
implementation of our argument.

When $\ell=k/2$ we will use the result from \cite{FK1} as our main
technical tool, and the above stated definition of $(\a,\e)$-regular
graphs. Here the definition of a pseudo-random hypergraph is
explicitly tailored to our application. Let $H=(V,\cE)$ be a
$k$-uniform hypergraph with vertex set $V=[n]$. Let
$\cP=(X_1,X_2,\ldots,X_{\n_\ell})$ be a partition of $[n]$ into
$\n_\ell$ parts each of size $\ell$. The graph $G_\cP=G_\cP(H)$ has
vertex set $[\n_\ell]$ and an edge $(i,j)$ whenever $E=X_i\cup
X_j\in \cE(H)$.
We now say that $H$ is {\em $(\a,\e)$-regular} if for a {\em
randomly chosen} $\cP$, the graph $G_\cP$ is $(\a,\e)$-regular \qs
\footnote{An event $\cE_n$ occurs {\em quite surely}, or \qs\ for
brevity, if $\Pr(\cE_n)=1-O(n^{-C})$ for any positive constant
$C$.}.
\begin{theorem}\label{th3}
Let $H=([n],\cE)$ be a $(p,\e)$-regular $k$-uniform hypergraph with
$k=2\ell$ and
$$\e^4np\gg \log^2n\ and\ \e^5np\gg\log(1/\e)\log n.$$
Then  $H$ contains a collection of $(1-20\e)\binom{n}{k}p/\n_\ell$
edge disjoint type $\ell$ Hamilton cycles.
\end{theorem}

We finally consider the case $k=\ell$. Here we will be packing
perfect matchings as opposed to Hamilton cycles. Let
$k_X=\rdown{k/2}$ and $k_Y=\rdup{k/2}$.
\begin{description}
\item[$R_a$:] $\displaystyle{\min_{S\in \binom{[n]}{k_X}}d_H(S)\geq\rdup{(1-\e)\binom{n-k_X}{k_Y}p}}$.
\item[$R_b$:] $\displaystyle{\min_{S\in \binom{[n]}{k_Y}}\,d_H(S)\geq\rdup{(1-\e)\binom{n-k_Y}{k_X}p}}$.
\item[$R_c$:] $\displaystyle{\max_{S\in \binom{[n]}{k_X+1}}\,d_H(S)\leq\rdup{(1+\e)\binom{n-k_X-1}{k_Y-1}p}}$.
\item[$R_d$:] $\displaystyle{\max_{S\in \binom{[n]}{k_Y+1}}\,d_H(S)\leq\rdup{(1+\e)\binom{n-k_Y-1}{k_X-1}p}}$.
\item[$R_e$:] $\displaystyle{\max_{\substack{S_1,S_2\in \binom{[n]}{k_X}\\ S_1\cap S_2=\emptyset}}}\,|N_H(S_1)\cap N_H(S_2)|
\leq \rdown{(1+\e)\binom{n-2k_X}{k_Y}p^2}$.
\item[$R_f$:] $\displaystyle{\max_{\substack{S_1,S_2\in \binom{[n]}{k_Y}\\ S_1\cap S_2=\emptyset}}}\,|N_H(S_1)\cap N_H(S_2)|
\leq \rdown{(1+\e)\binom{n-2k_Y}{k_X}p^2}$.
\end{description}

\begin{theorem}\label{th4}
Let $H=([n],\cE)$ be a $k$-uniform hypergraph with $m$ edges that
satisfies $\cR=\{R_a,R_b,R_c,R_d,R_e,R_f\}$ and suppose that $1\gg \e\gg
\log^5n/(n^{1/2}p^2)$. Then  $H$ contains a collection of
$(1-4\e^{1/3})m/\n_k$ edge disjoint perfect matchings.
\end{theorem}

Am interesting point of reference for our theorems is results about
{\em perfect} decompositions of the edge set of a complete
$k$-uniform hypergraph $K_n^k$ into Hamilton cycles of various types
(assuming of course some natural divisibility conditions). These
include a recent result of Bailey and Stevens \cite{BS} about
packing tight Hamilton cycles and a famous result of Baranyai
\cite{Bar} about decomposing the edge set of $K_n^k$ into perfect
matchings. While we do not -- and can not for obvious reasons --
achieve perfect decomposition, but rather pack almost all edges, our
results apply to a wide class of hypergraphs, including relatively
sparse hypergraphs.

In the next section we focus on $H=H(n,p,k)$ and first prove Theorem
\ref{th5} for $k=3$. We then give a proof for general $k$. In
Section \ref{secth1} we prove Theorems \ref{th1}, \ref{th3} and
\ref{th4}. The last section is devoted to concluding remarks.

\section{Random hypergraphs}\label{hnp3}
We prove Theorem \ref{th5} in this section.

The proof for $2\ell>k$ is based on the same idea as for the case
$k=3,\ell=2$ but is heavier on notation and will be given
immediately afterwards. Hopefully, the reader will find it useful to
consider the simplest case first. The proof for random hypergraphs
is simpler than the proof for regular (i.e., pseudo-random)
hypergraphs and hopefully will help in the understanding of the
proofs of Theorems \ref{th1} and \ref{th3}.

{\bf Case 1:} $k=3,\ell=2$.\\
We will construct the Hamilton cycles via the following algorithm:
\begin{description}
\item[$A_1$:] Choose $r=n(np)^{1/2}$ random partitions
$(X_i,Y_i),i=1,2,\ldots r$, of $V$ into two sets of size $\n_2$.

We use the notation
$$X_i=\set{x_{i,1}<x_{i,2}<\cdots<x_{i,\n_2}}\ and\ Y_i=\set{y_{i,1}<y_{i,2}<\cdots<y_{i,\n_2}}.$$
For each $i$ we choose a random permutation $\s_i$ on $X_i$ and define a Hamilton cycle
$$\G_{i}=(x_{i,\s_i(1)},x_{i,\s_i(2)},\ldots,x_{i,\s_i(\n_2)}).$$
\item[$A_2$:] At this point we expose the edges of $H(n,p,3)=([n],\cE_p)$.
\item[$A_3$:] Suppose now that for edge $E\in \cE_p$ there are $f(E)$ instances
$i$ such that $|X_i\cap E|=2$ and $X_i\cap E$ is an edge of
$\G_{i}$. If $f(E)>0$, then choose one of the $f(E)$ instances at
random and label the edge $E$ with the chosen $i$; if $f(E)=0$, the
edge $E$ stays unlabelled. Let $H_i\subset H$ be the subhypergraph
of all edges labelled by $i$.
\item[$A_4$:] Let $G_{i}$ be the bipartite graph with
vertex set $A_{i}\cup Y_i$ defined as follows: $A_{i}$ is a copy of
$[\n_2]$ (viewed as the set of edges of the Hamilton cycle $\G_i$).
Add edge $(a,b)$ to $G_{i}$ if
$E=(x_{i,\s_i(a)},y_{i,b},x_{i,\s_i(a+1)})\in H_i$ (i.e., $(a,b)\in
G_i$ if the $a$-th edge of the cycle $\G_i$ united with the vertex
$y_b$ forms an edge $E$ of $H$ labeled by $i$).
\item[$A_5$:] We claim that \whp\ (see Lemma \ref{lem1} below)
each $G_{i}$ will contain at least
$$n_0=(1-\e)\n_2 p_0$$
edge disjoint
perfect matchings.

Here
\begin{equation}\label{k0}
\e=\bfrac{72n\log n}{r}^{1/2}\text{ and } f_0=
r\r+(12r\r\log n)^{1/2}\text{ and
}p_0=\frac{p}{f_0}.
\end{equation}
Here $\r=\r_{3,2}$ where
$$\r_{k,\ell}=\frac{\n_\ell^2}{\binom{n}{k}}.$$
\begin{remark}\label{rem2}
Note that $\r$ is the probability that instance $i$ is one of the $f(E)$ instances in
$A_3$ for edge $E\in \cE$. 
\end{remark}

Each such matching gives rise to a loose Hamilton cycle of $H_i$ and
these will be edge disjoint by construction. Indeed suppose that our
matching is $(e_a,\f(e_a)),a=1,2,\ldots,\n_2$, where the edges $e_a$
are ordered according to the order of their appearance along the
Hamilton cycle $\G_i$. From this we obtain the type 2 Hamilton cycle
with edges $E_a=e_a\cup \set{\f(e_a)}$. Since the subhypergraphs
$H_i$ are edge disjoint and since distinct edges in the graph $G_i$
correspond to distinct edges of $H_i$, the so obtained Hamilton
cycles in $H$ are indeed edge disjoint.
\end{description}
It follows that \whp\ $H(n,p,3)$ contains at least $rn_0$ edge
disjoint Hamilton cycles, proving Theorem \ref{th5} for this case.

\begin{lemma}\label{lem1}
$$\Pr(G_{i}\text{ does not contain }n_0\text{ edge disjoint perfect matchings})=o(n^{-3}) .$$
\end{lemma}
\proofstart
The Max-Flow Min-Cut theorem
tells us that the following is a necessary and sufficient condition for $\G_{i}$ to have $n_0$
edge disjoint perfect matchings: Suppose that we make up a network with source $\s$ and sink $\t$
and join $\s$ to each vertex of $A=A_{i}$ by an edge of capacity $n_0$ and each vertex of $B=Y_i$ to $\t$ by an edge of capacity
$n_0$. Each edge of $G=\G_{i}$ is given capacity one. Suppose that our minimum cut is $X:\bar{X}$ and $S=A\cap X$ and $T=B\cap X$
then a necessary and sufficient condition for the existence of $n_0$ disjoint perfect matchings is that
$$(\n_2-|S|)n_0+|T|n_0+e(S,B\setminus T)\geq n_0\n_2$$
which reduces to
\begin{equation}\label{1}
m\geq (s-t)n_0
\end{equation}
for all $S\subseteq A,T\subseteq B$, if $|S|=s,|T|=t,m=e(S,B\setminus T)$.

Note that we need only verify \eqref{1} computationally for $t\leq
\n_2/2$. When $t>\n_2/2$ we could repeat our computations to show
that \whp\ $e(B\setminus T,A\setminus (A\setminus S))\geq
((\n_2-t)-(\n_2-s))n_0$.

 For a triple $E\subset [n]$, we say that $1\le i\le r$ {\em includes}
 $E$ if the set $E\cap X_i$ is of size 2 and is one of the edges
 of the cycle $\G_i$. Thus the random variable $f(E)$ counts the
 number of partitions $(X_i,Y_i)$ that include $E$. Observe that the
 $i$-th partition includes a fixed triple $E$ with probability
 $$
 \binom{3}{2}\frac{\binom{n-3}{\n_2-2}}{\binom{n}{\n_2}}\,\frac{\n_2}{\binom{\n_2}{2}}
 =\frac{3n}{2(n-1)(n-2)}=\frac{\n_2^2}{\binom{n}{3}}
$$
(first choose two elements of $E\cap X_i$, then choose $X_i$ to
intersect $E$ in exactly these two elements, then choose a Hamilton
cycle in $X_i$ -- due to symmetry the probability that $E\cap X_i$
is one of its $\n_2$ edges is $\n_2\binom{\n_2}{2}^{-1}$). Moreover
the events ``$i$ includes $E$" are mutually independent for
different $i$. Therefore, the random variable $f(E)$ is distributed
binomially with parameters $r$ and $\r$. Now
using the following Chernoff bounds for $0\leq\e\leq1$:
\begin{eqnarray}
\Pr(Bin(m,\xi)-m\xi\leq -\e m\xi)&\leq& e^{-\e^2m \xi/2}\label{chern1}\\
\Pr(Bin(m,\xi)-m\xi\geq \e m\xi)&\leq& e^{-\e^2m
\xi/3}\label{chern2}
\end{eqnarray}
we see that with probability at least $1-o(1)$ we have $1\leq f(E)\leq
f_0$ for all $\binom{n}{3}$ possible edges. So assume that indeed
$1\leq f(E)\le f_0$ for all $E$. Moreover, the values of $f(E)$ are
determined by Steps $A_1$ and $A_2$ of our construction and are thus
independent of the appearance of random edges at Step $A_3$. For
$1\le a,b \le \n_2$, the pair $(a,b)$ is an edge of the random
auxiliary graph $G_i$ if the corresponding triple $E$ is an edge of
the random hypergraph $H$ and is chosen to be labelled by $i$. Thus
$(a,b)\in E(G_i)$ independently and with probability at least
$p/f_0=p_0$.

Therefore we can \whp\ reduce our problem to showing that \whp\ the random bipartite graph $K_{\n_2,\n_2,p_0}$
contains $n_0$ edge disjoint perfect matchings.

Then with $\e$ as defined in \eqref{k0},
\begin{align}
&Pr(\exists S\subseteq A,T\subseteq B,|S|> |T|,|T|\leq \n_2/2:\;e(S,B\setminus T)\leq (1-\e)|S|(\n_2-|T|)p_0)\leq\nonumber\\
&\sum_{s=1}^{\n_2}\sum_{t=1}^{\min\set{s-1,\n_2/2}}\binom{\n_2}{s}\binom{\n_2}{t}
\exp\set{-\frac{\e^2}{2}s(\n_2-t)p_0}\leq\nonumber\\
&\sum_{s=1}^{\n_2}\sum_{t=1}^{\min\set{s-1,\n_2/2}}\binom{\n_2}{s}\binom{\n_2}{t}\exp\set{-\e^2s\n_2p_0/4}.\label{sum}
\end{align}
Assume first that $\binom{\n_2}{s}\geq \binom{\n_2}{t}$. Then
$$\eqref{sum}\leq \sum_{s=1}^{\n_2}\sum_{t=1}^{\min\set{s-1,\n_2/2}}\brac{\frac{n^2e^2}{4s^2}\cdot e^{-\e^2\n_2p_0/4}}^s
=\sum_{s=1}^{\n_2}\sum_{t=1}^{\min\set{s-1,\n_2/2}}\brac{\frac{n^2e^2}{4s^2}\cdot n^{-6+o(1)}}^s=o(n^{-3}).$$
When $\binom{\n_2}{s}\leq \binom{\n_2}{t}$ we can replace \eqref{sum} by
$$ \sum_{s=\n_2/2}^{\n_2}\sum_{t=1}^{\min\set{s-1,\n_2/2}}\brac{\frac{n^2e^2}{4t^2}\cdot n^{-6+o(1)}}^t=o(n^{-3}).$$
Here we have used $s\geq t$.
\proofend

It follows (see \eqref{1}) that \whp\
$$m\geq (1-\e)s(\n_2-t)np_0\geq (1-\e)(s-t)np_0/2=(s-t)nn_0.$$
This completes the proof of Theorem \ref{th1} for $k=3,\ell=2$.

\bigskip

With a roadmap in mind, we proceed to the general case.

{\bf Case 2:} $\ell<k<2\ell$.\\
We will construct the Hamilton cycles via the following algorithm:
\begin{description}
\item[$B_1$:] Let
$$N_X=\frac{k-\ell}{\ell}n\ and\ N_Y=\frac{2\ell-k}{\ell}n.$$
Choose $r=n^{k-2}(np)^{1/2}$ random partitions
$(X_i,Y_i),i=1,2,\ldots r$, of $[n]$ into two sets of size $N_X$ and
$N_Y$ respectively.

We use the notation
$$X_i=\set{x_{i,1}<x_{i,2}<\cdots<x_{i,N_X}}\ and\ Y_i=\set{y_{i,1}<y_{i,2}<\cdots<y_{i,N_Y}}.$$
\item[$B_2$:] At this point we expose the edges of $H(n,p,k)=([n],\cE_p)$.
\item[$B_3$:] For each $i$ we let $\s_i$ be a random permutation of $X_i$ and let
$\t_i$ be a random permutation of $Y_i$. Form the partition
$X_{i,a},a=1,2,\ldots,\n_\ell$, of $X_i$ into sets of size $k-\ell$
and the partition $Y_{i,b},b=1,2,\ldots,\n_\ell$, of $Y_i$ into sets
of size $2\ell-k$. Here
$X_{i,a}=\set{x_{i,\s_i((a-1)(k-\ell)+1)},\ldots,x_{i,\s_i(a(k-\ell))}}$
and
$Y_{i,b}=\set{y_{i,\t_i((b-1)(2\ell-k)+1)},\ldots,y_{i,\t_i(b(2\ell-k))}}$.

We define the ``Hamilton cycle''
$$\G_{i}=(X_{i,1},X_{i,2}\ldots,X_{i,\n_\ell}).$$

\item[$B_4$:] Suppose now that for $E\in \cE_p$ there are $f(E)$ instances
$i$ such that for some $a,b$ and some partition $S_1,S_2,S_3$ of $E$
we have $S_1=X_{i,a},S_2=X_{i,a+1}$ (where we set $\n_\ell+1$ to be
equal to $1$) and $S_3=Y_{i,b}$. We say that $i$ includes $E$.
Choose one of the $f(E)$ instances at random and label edge $E$ with
the chosen $i$. If $f(E)=0$, the edge $E$ stays unlabeled. Let $H_i$
be the subhypergraph of $H$ formed by the edges of $H$ labeled by
$i$.

\item[$B_5$:] Let $G_{i}$ be the bipartite graph with vertex partition $A_{i}$ and
$B_i$ comprising disjoint copies of $[\n_\ell]$. For $a\in A_{i}$
and $b\in B_i$  we make $(a,b)$ an edge of $G_{i}$ if $E=X_{i,a}\cup
Y_{i,b}\cup X_{i,a+1}\in \cE_p$ and $E$ is labelled with $i$.
\item[$B_6$:] We claim (see Lemma \ref{lem1x}) that \whp\
each $G_{i}$ will contain at least
$$n_0=(1-\e)\n_\ell p_0$$
edge disjoint
perfect matchings.

Here
\begin{equation}\label{k0x}
\e=\bfrac{4(k+3)k!\log n}{\ell(np)^{1/2}}^{1/2}\text{ and }f_0= \r r+(4k\r
r\log n)^{1/2}\text{ and }p_0=\frac{p}{f_0},
\end{equation}
where $\r=\r_{k,\ell}$.
\end{description}

Each such matching
gives rise to a type $\ell$ Hamilton cycle of $H$ and these will be edge disjoint by
construction.

In this way we obtain at least
$rn_0$
edge disjoint Hamilton cycles, proving Theorem \ref{th5} for the case $\ell<k<2\ell$.
\proofend

\begin{lemma}\label{lem1x}
$$\Pr(G_{i}\text{ does not contain }n_0\text{ edge disjoint perfect matchings})=o(n^{-k}) .$$
\end{lemma}
\proofstart The edges of $G_{i}$ appear independently with
probability $p/f(E)$ where $f(E)$ has distribution $Bin(r,\r)$. (To
see it, for a fixed partition $(X_i,Y_i)$ and a fixed pair of
permutations $(\s_i,\t_i)$ of $X_i,Y_i)$, resp., the index $i$
includes $\n_\ell^2$ $k$-tuples from $[n]$. Therefore by symmetry a
random $i$ includes a fixed $k$-tuple $E$ with probability
$\frac{\n_\ell^2}{\binom{n}{k}}=\r_{k,\ell}$.) So \eqref{chern1} and \eqref{chern2} imply that $1\leq f(E)\leq f_0$ 
with probability $1-O(n^{-4k/3})$.

We have reduced our problem to showing that \whp\ the random
bipartite graph $K_{\n_\ell,\n_\ell,p_0}$ contains $n_0$ edge
disjoint perfect matchings. We need to verify \eqref{1}
computationally for $t\leq \n_\ell/2$. Then with $\e$ as defined in
\eqref{k0x},
\begin{align}
&Pr(\exists S\subseteq A,T\subseteq B,|S|\geq |T|,|T|\leq \n_\ell/2:\;e(S,B\setminus T)\leq (1-\e)|S|(\n_\ell-|T|)p_0)\leq\nonumber\\
&\sum_{s=1}^{\n_\ell}\sum_{t=1}^{\min\set{s-1,\n_\ell/2}}\binom{\n_\ell}{s}
\binom{\n_\ell}{t}\exp\set{-\frac{\e^2}{2}s(\n_\ell-t)p_0}\leq\nonumber\\
&\sum_{s=1}^{\n_\ell}\sum_{t=1}^{\min\set{s-1,\n_\ell/2}}\binom{\n_\ell}{s}\binom{\n_\ell}{t}\exp\set{-\e^2s\n_\ell p_0/4}.\label{sumx}
\end{align}
Assume first that $\binom{\n_\ell}{s}\geq \binom{\n_\ell}{t}$. Then
$$\eqref{sumx}\leq \sum_{s=1}^{\n_\ell}\sum_{t=1}^{\min\set{s-1,\n_\ell/2}}\brac{\frac{n^2e^2}{\ell^2s^2}\cdot
e^{-\e^2\n_\ell p_0/4}}^s
\le\sum_{s=1}^{\n_\ell}\sum_{t=1}^{\min\set{s-1,\n_\ell/2}}\brac{\frac{n^2e^2}{\ell^2s^2}\cdot
n^{-k-3}}^s=o(n^{-k}).$$ When $\binom{\n_\ell}{s}\leq
\binom{\n_\ell}{t}$ we can replace \eqref{sumx} by
$$ \sum_{s=\n_\ell/2}^{\n_\ell}\sum_{t=1}^{\min\set{s-1,\n_\ell/2}}\brac{\frac{n^2e^2}{\ell^2t^2}\cdot
n^{-k-3}}^t=o(n^{-k}).$$ Here we have used $s\geq t$. \proofend

{\bf Case 3:}  $k=2\ell$.

When $\ell=k/2$ we have $N_Y=0$ and the argument above breaks down. We can however use our result from \cite{FK1}
to obtain something.

We will construct the Hamilton cycles via the following algorithm:
\begin{description}
\item[$C_1$:]
Choose $r=n^{k-2}(np)^{1/2}$ random partitions $\cP_i,i=1,2,\ldots
r$, of $[n]$ into $\n_\ell$ sets $X_{i,a}$ of size $\ell$.
\item[$C_2$:] Expose the edges of $H(n,p,k)=([n],\cE_p)$.
\item[$C_3$:] For each $E\in \cE_p$ we let $f(E)$ denote the number of partitions
$i$ such that $\cP_i$ contains a pair of parts $X_{i,a},X_{i,b}$
such that $X_{i,a}\cup X_{i,b}=E$. The random variable $f(E)$ is
distributed as $Bin(r,\r)$ where
$$\r=\frac{\binom{\n_\ell}{2}}{\binom{n}{k}}.$$
So \eqref{chern1} and \eqref{chern2} imply that 
$(1-\e)f_0\le f(E)\leq f_0$ with probabilty $1-O(n^{-4k/3})$. Choose one of these $f(E)$ instances
at random and label the edge $E$ with the chosen $i$.  Let $H_i$ be
the subhypergraph of all edges of $H$ labeled by $i$. Here we can
use $\e,f_0,p_0$ as in \eqref{k0x}.
\item[$C_5$:] For each $i$ let $G_{i}$ be the graph with vertex set
$[\n_\ell]$, where $a,b\in [\n_\ell]$ are connected by an edge if
$X_{i,a}\cup X_{i,b}$ is an edge of $H$ labeled by $i$. We will show
below in Lemma \ref{regux} that \qs\ each $G_i$ is
$\brac{(1-2\e)p_0,2\e p_0}$-regular.
\item[$C_6$:] We then apply Theorem \ref{FK1} to show that \qs\
each $G_i$ contains at least $\brac{(1-2\e)p_0/2-8\e p_0}\n_\ell$
edge disjoint Hamilton cycles. Each such Hamilton cycle corresponds
to a Hamilton cycle of type $\ell$ in $H_i$, and the
Hamilton cycles so obtained are edge disjoint.
\end{description}
Thus $H$ contains at least
$$r\brac{(1-2\e)p_0/2-8\e p_0}\n_\ell\geq
\binom{n}{k}\frac{p}{\n_\ell}\brac{1-20\e}$$ edge disjoint type
$\ell$ Hamilton cycles, completing the proof of Theorem \ref{th5}
for this case.

\begin{lemma}\label{regux}
Each $G_i$ is \qs\ $\brac{(1-2\e)p_0,2\e p_0}$-regular.
\end{lemma}
\proofstart The degree of vertex $v$ in $G_i$ dominates
$Bin(\n_\ell-1,p_0)$ and so Property $Q_a$ holds from Chernoff
bounds. Observe that $\n_\ell p_0=\Omega((np)^{1/2})\gg\log n$.
Similarly the number of edges between two sets $S,T$ dominates
$Bin(|S|\,|T|,p_0)$ and is dominated by $Bin(|S|\,|T|,(1-\e)^{-1}p_0)$
and Property $Q_b$ also holds from Chernoff bounds. \proofend

{\bf Case 4:}  $k=\ell$.

Here the aim is to find many edge disjoint perfect machings. We construct them via the following algorithm:
\begin{description}
\item[$D_1$:] Let $k_X=\rdown{k/2}$ and $k_Y=\rdup{k/2}$ and
$$N_X=\frac{k_X}{k}n\ and\ N_Y=\frac{k_Y}{k}n.$$
\item[$D_2$:]
Choose $r=n^{k-2}(np)^{1/2}$ random partitions
$(X_i,Y_i),i=1,2,\ldots r$, of $[n]$ into two sets of size $N_X$ and
$N_Y$ respectively.

We use the notation
$$X_i=\set{x_{i,1}<x_{i,2}<\cdots<x_{i,N_X}}\ and\ Y_i=\set{y_{i,1}<y_{i,2}<\cdots<y_{i,N_Y}}.$$
\item[$D_3$:] At this point we expose the edges of $H(n,p,k)=([n],\cE_p)$.
\item[$D_4$:] For each $i$ we let $\s_i$ be a random permutation of $X_i$ and let
$\t_i$ be a random permutation of $Y_i$. Form the partition
$X_{i,a},a=1,2,\ldots,\n_k$, of $X_i$, into sets of size $k_X$ and
the partition $Y_{i,b},b=1,2,\ldots,\n_k$, of $Y_i$ into sets of
size $k_Y$. Here
$X_{i,a}=\set{x_{i,\s_i((a-1)k_X+1)},\ldots,x_{i,\s_i(ak_X)}}$ and
$Y_{i,b}=\set{y_{i,\t_i((b-1)k_Y+1)},\ldots,y_{i,\t_i(bk_Y)}}$.

\item[$D_5$:] Suppose now that for $E\in \cE_p$ there are $f(E)$ instances
$i$ such that for some $a,b$ and some partition $S_1,S_2$ of $E$ we
have $S_1=X_{i,a}$ and $S_2=Y_{i,b}$. We say that $i$ includes $E$.
Choose one of the $f(E)$ instances at random and label edge $E$ with
the chosen $i$.
\item[$D_6$:] Let $G_{i}$ be the bipartite graph with vertex partition $(A_{i},B_i)$
 comprising disjoint copies of $[\n_k]$. For $a\in A_{i}$ and
$b\in B_i$  we make $(a,b)$ an edge of $G_{i}$ if $E=X_{i,a}\cup
Y_{i,b}\in \cE_p$ and $E$ is labelled with $i$. So, by construction,
each $E\in \cE_p$ is associated with at most one $G_{i}$.
\item[$D_7$:] We claim (see Lemma \ref{lem1y}) that \whp\
each $G_{i}$ will contain at least
$$n_0=(1-\e)\n_k p_0$$
edge disjoint
perfect matchings.

Here
\begin{equation}\label{k0y}
\e=10k!\bfrac{4(k+3)k!\log n}{\ell(np)^{1/2}}^{1/2}\text{ and }f_0= \r
r+(4k\r r\log n)^{1/2}\text{ and }p_0=\frac{p}{f_0},
\end{equation}
where $\r=\r_{k,k}$.
\end{description}
Thus \whp\ $H$ contains at least $rn_0$ edge disjoint perfect matchings and this completes the proof of Theorem \ref{th5}.
\begin{lemma}\label{lem1y}
$$\Pr(G_{i}\text{ does not contain }n_0\text{ edge disjoint perfect matchings})=o(n^{-k}) .$$
\end{lemma}
\proofstart The edges of $G_{i}$ appear independently with
probability $p/f(E)$ where $f(E)$ has distribution $Bin(r,\r)$. So \eqref{chern1} and \eqref{chern2} imply
that  $1\leq f(E)\leq
f_0$ with probability $1-O(n^{-4k/3})$.

We have reduced our problem to showing that \whp\ the random
bipartite graph $K_{\n_k,\n_k,p_0}$ contains $n_0$ edge disjoint
perfect matchings. We need to verify \eqref{1} computationally for
$t\leq \n_k/2$. We follow the proof of Lemma \ref{lem1x} with
$\ell=k$. \proofend

\section{Pseudo-random hypergraphs}\label{secth1}
In this section we prove Theorems \ref{th1}, \ref{th3} and
\ref{th4}. We follow the same strategy as described in Section
\ref{hnp3}. There are complications caused by the notation that we
have to add and also by the fact that $H$ is not random.

{\bf Case 1:} $\ell<k<2\ell$ (Theorem \ref{th1}).

We will construct the Hamilton cycles via the following algorithm: First choose $f_0$ such that
\beq{f0l}
\frac{\log^2n}{\e^4}\ll f_0^2\ll \frac{\e n^{1/2}p^2}{\log n}.
\eeq
\begin{description}
\item[$E_1$:] Let
$$r=(1-\e)\binom{n}{k}\frac{f_0}{\n_\ell^2};\ \ p_0=\frac{p}{f_0}\,.$$
\item[$E_2$:]
Let
$$N_X=\frac{k-\ell}{\ell}n\ and\ N_Y=\frac{2\ell-k}{\ell}n.$$
Now choose $r$ random partitions $(X_i,Y_i),i=1,2,\ldots r$, of
$[n]$ into two sets of size $N_X$ and $N_Y$ respectively.

We use the notation
$$X_i=\set{x_{i,1}<x_{i,2}<\cdots<x_{i,N_X}}\ and\ Y_i=\set{y_{i,1}<y_{i,2}<\cdots<y_{i,N_Y}}.$$
\item[$E_3$:] For each $i$ we let $\s_i$ be a random permutation of $X_i$ and let
$\t_i$ be a random permutation of $Y_i$. Form the partition
$X_{i,a},a=1,2,\ldots,\n_\ell$, of $X_i$ into sets of size $k-\ell$
and the partition $Y_{i,b},b=1,2,\ldots,\n_\ell$, of $Y_i$ into sets
of size $2\ell-k$. Here
$X_{i,a}=\set{x_{i,\s_i((a-1)(k-\ell)+1)},\ldots,x_{i,\s_i(a(k-\ell))}}$
and
$Y_{i,b}=\set{y_{i,\t_i((b-1)(2\ell-k)+1)},\ldots,y_{i,\t_i(b(2\ell-k))}}$.

We define the ``Hamilton cycle''
$$\G_{i}=(X_{i,1},X_{i,2},\ldots,X_{i,\n_\ell}).$$
\item[$E_4$:] Suppose now that for $E\in \cE$ there are $f(E)$ instances
$i$ such that for some $a,b$ and some partition $S_1,S_2,S_3$ of $E$
we have $S_1=X_{i,a},S_2=X_{i,a+1}$ and $S_3=Y_{i,b}$.
 Choose one of the $f(E)$ instances at random
and label edge $E$ with the chosen $i$.

Thus $f(E)$ is distributed as $Bin(r,\r)$ where $\r=\r_{k,\ell}$.
So \eqref{chern1} and \eqref{chern2} imply that\\
$\big(1-\frac{3}{2}\e\big)f_0\leq f(E)\leq f_0$ with probability $1-o(n^{-k})$.
\item[$E_5$:] Let $G_{i}$ be the bipartite graph with vertex partition $A_{i}$ and
$B_i$ comprising disjoint copies of $[\n_\ell]$. For $a\in A_{i}$
and $b\in B_i$  we make $(a,b)$ an edge of $G_{i}$ if $E=X_{i,a}\cup
Y_{i,b}\cup X_{i,a+1}\in \cE$ and $E$ is labelled with $i$.
\item[$E_6$:] We claim that \whp\ (see Lemma \ref{lem2} below)
each $G_{i}$ will contain at least
$$n_0=(1-(5\e)^{1/3})\n_\ell p_0$$
edge disjoint
perfect matchings.
\end{description}

Each such matching
gives rise to a type $\ell$ Hamilton cycle of $H$ and these will be edge disjoint by
construction.

In this way we obtain at least
$$rn_0\geq \binom{n}{k}\frac{p}{\n_\ell}(1-(5\e)^{1/3}-\e)$$
edge disjoint Hamilton cycles, proving Theorem \ref{th1} for the case $\ell<k<2\ell$.

We will show later (Lemma \ref{lem2} below) that if we can prove
that the degrees and co-degrees of our bipartite graphs $G_{i}$
``behave'', then we can deduce the existence of many disjoint
perfect matchings and so get our packing of Hamilton cycles. Given
Lemma \ref{lem2}, all we need to do is to estimate the degrees and
co-degrees of vertices in a fixed $G_{i}$.
\begin{lemma}\label{lem3}
{\bf Whp}, over our random choices of $X_i,Y_i,\s_i,\t_i$, each
$G_{i}$ has minimum degree at least $(1-2\e)\n_\ell p_0$ and maximum
co-degree of at most $(1+5\e)\n_\ell p_0^2$.
\end{lemma}
\proofstart We fix $i$  and focus on $G_{i}$. We first show that the
minimum degree in $G_{i}$ is large. We first fix $a\in A_{i}$. The
vertex $a$ corresponds to the block $X_{i,a}$ of $\sigma_i$.
Condition on $X_{i,a}\cup X_{i,a+1}=S$ for some $S\subset [n]$,
$|S|=2(k-\ell)$. We expose a random subset $Y_i$ first. Let $Z_a^*$
be the number of edges $E\in \cE$ such that $S\subset E$ and $E\cap
Y_i=E-S$. For each edge $E\in N_H(S)$
$$
\Pr(E\cap Y_i=E-S)
=\frac{\binom{n-k}{N_Y-(2\ell-k)}}{\binom{n-2(k-\ell)}{N_Y}}=\frac{(N_Y)_{2\ell-k}}{(n-2(k-\ell))_{2\ell-k}}=
\brac{1+O\bfrac{1}{n}}\bfrac{2\ell-k}{\ell}^{2\ell-k}\ .
$$
Therefore by assumption $P_a$,
$$
\E(Z_a^*)\ge
(1-\e)\brac{1+O\bfrac{1}{n}}\bfrac{2\ell-k}{\ell}^{2\ell-k}
\binom{n-2(k-\ell)}{2\ell-k}p\ .
$$
Since changing the fate of one vertex with respect to the choice of
$Y_i$ changes the value of $Z_a^*$ by at most
$$
\D_a=\max_{S'\in \binom{[n]}{2(k-\ell)+1}}d_{H}(S')\,,
$$
and the latter quantity is bounded by $(1+\e)\binom{n}{2\ell-k-1}p$
by assumption $P_c$, we get by the Azuma-Hoeffding inequality that
for any $t>0$
\beq{Zastar}
\Pr(Z_a^*\le \E(Z_a^*)-t)\le \exp
\left\{-\frac{2t^2}{n\D_a^2}\right\}\ .
\eeq
Here we are using the following inequality: Let $S_n$ denote the set of permutations of $[n]$ and let
$f:S_n\to\Re$ be such that $|f(\p)-f(\p')|\leq u$ whenever $\p'$ is obtained from $\p$ by transposing 
two elements. Then if $\p$ is chosen randomly from $S_n$ then
\beq{ranpi}
\Pr(f(\p)-\E(f)\leq -t)\leq \exp\set{-\frac{2t^2}{nu^2}}.
\eeq
For a proof see e.g., Section 3.2 of
\cite{McD} or Lemma 11 of \cite{FP}.

In this context, think of choosing a random $m$-subset of $[n]$ as chosing a random $\p$
and then taking the first $m$ elements as your subset.

Plugging in the estimates on $\E(Z_a^*)$ and $\D_a$ stated above in \eqref{ranpi}, we
get that \qs\ for every $a\in A_i$,
\begin{eqnarray}
Z_a^* &\ge&
(1-\e)\brac{1+O\bfrac{1}{n}}\bfrac{2\ell-k}{\ell}^{2\ell-k}
\binom{n-2(k-\ell)}{2\ell-k}p-n^{2\ell-k-1/2}p\log n\nonumber\\
&\ge&
\big(1-\frac{3}{2}\e\big)\bfrac{2\ell-k}{\ell}^{2\ell-k}\binom{n}{2\ell-k}p\
.\label{Za*}
\end{eqnarray}

So assume that $Y_i$ is chosen so that (\ref{Za*}) holds. Now we
expose the random permutation $\t_i$ of $Y_i$. Let $Z_a$ be the
degree of $a$ in $G_i$, which is the number of edges $E\in\cE$ such
that
\begin{enumerate}
\item $E\cap Y_i=E-S$ (the number of such edges is $Z_a^*$);
\item $E\cap Y_i$ forms a block $Y_{i,b}$ under $\t_i$;
\item $E$ is labeled by $i$ (this happens independently and with
probability $1/f(E)\ge 1/f_0$).
\end{enumerate}
Hence,
$$
\E(Z_a)\ge Z_a^*\, \frac{\n_{\ell}}{\binom{N_Y}{2\ell-k}}\,
\frac{1}{f_0}\ .
$$
Observe that changing $\t_i$ by a single transposition changes the
value of $Z_a$ by at most 2 (at most two blocks $Y_{i,b}$ are
affected by such a change). Therefore, applying concentration
results for permutation graphs we get that for any $t>0$
$$
\Pr(Z_a\le \E(Z_a)-t)\le \exp\left\{-\frac{t^2}{2N_Y}\right\}\ .
$$
Thus \qs\ for every partition $i$ and for every $a\in A_i$, its
degree is $G_i$ is at least
$$
\big(1-\frac{3}{2}\e\big)\bfrac{2\ell-k}{\ell}^{2\ell-k}\binom{n}{2\ell-k}p\,
\frac{\n_{\ell}}{\binom{N_Y}{2\ell-k}}\, \frac{1}{f_0} -n^{1/2}\log
n \ge (1-2\e)\n_{\ell}p_0\,,
$$
due to our assumption on $\e$.

The argument for the degrees of the vertices of
$B_i$ is quite similar.  Fix $b\in B_{i}$. The vertex $b$
corresponds to the block $Y_{i,b}$ of $\t_i$. Condition on
$Y_{i,b}=S$ for some $S\subset [n]$, $|S|=2\ell-k$. We expose a
random subset $X_i$ first. Let $Z_b^*$ be the number of edges $E\in
\cE$ such that $S\subset E$ and $E\cap X_i=E-S$. For each edge $E\in
N_H(S)$
$$
\Pr(E\cap X_i=E-S)
=\frac{\binom{n-k}{N_X-2(k-\ell)}}{\binom{n-(2\ell-k))}{N_X}}=\frac{(N_X)_{2k-2\ell}}{(n-(2\ell-k))_{2k-2\ell}}=
\brac{1+O\bfrac{1}{n}}\bfrac{k-\ell}{\ell}^{2k-2\ell}\ .
$$
Therefore by assumption $P_b$,
$$
\E(Z_b^*)\ge
(1-\e)\brac{1+O\bfrac{1}{n}}\bfrac{k-\ell}{\ell}^{2k-2\ell}
\binom{n-2\ell+k}{2(k-\ell)}p\ .
$$
Since changing the fate of one vertex with respect to the choice of
$X_i$ changes the value of $Z_b^*$ by at most
$$
\D_b=\max_{S'\in \binom{[n]}{2\ell-k+1}}d_{H}(S')\,,
$$
and the latter quantity is bounded by $(1+\e)\binom{n}{2k-2\ell-1}p$
by assumption $P_d$, we get by \eqref{ranpi} that
for any $t>0$
$$
\Pr(Z_b^*\le \E(Z_b^*)-t)\le \exp
\left\{-\frac{2t^2}{n\D_b^2}\right\}\ .
$$
Plugging in the estimates on $\E(Z_b^*)$ and $\D_b$ stated above, we
get that \qs\ for every $b\in B_i$,
\begin{eqnarray}
Z_b^* &\ge&
(1-\e)\brac{1+O\bfrac{1}{n}}\bfrac{k-\ell}{\ell}^{2k-2\ell}
\binom{n-2\ell+k}{2(k-\ell)}p-n^{2k-2\ell-1/2}p\log n\nonumber\\
&\ge&
\big(1-\frac{3}{2}\e\big)\bfrac{k-\ell}{\ell}^{2k-2\ell}\binom{n}{2k-2\ell}p\
.\label{Zb*}
\end{eqnarray}

So assume that $X_i$ is chosen so that (\ref{Zb*}) holds. Now we
expose the random permutation $\s_i$ of $X_i$. Let $Z_b$ be the
degree of $b$ in $G_i$, which is the number of edges $E\in\cE$ such
that
\begin{enumerate}
\item $E\cap X_i=E-S$ (the number of such edges is $Z_b^*$);
\item $E\cap X_i$ forms  two consecutive blocks $X_{i,a}$, $X_{i,a+1}$ under $\s_i$;
\item $E$ is labeled by $i$ (this happens independently and with
probability $1/f(E)\ge 1/f_0$).
\end{enumerate}
Hence,
$$
\E(Z_b)\ge Z_b^*\, \frac{\n_{\ell}}{\binom{N_X}{2k-2\ell}}\,
\frac{1}{f_0}\ .
$$
Observe that changing $\s_i$ by a single transposition changes the
value of $Z_b$ by at most 4. Therefore, applying again concentration
results for permutation graphs  we get that for any $t>0$
$$
\Pr(Z_b\le \E(Z_b)-t)\le \exp\left\{-\frac{t^2}{8N_X}\right\}\ .
$$
Thus \qs\ for every partition $i$ and for every $b\in B_i$, its
degree is $G_i$ is at least
$$
\big(1-\frac{3}{2}\e\big)\bfrac{k-\ell}{\ell}^{2k-2\ell}\binom{n}{2k-2\ell}p\,
\frac{\n_{\ell}}{\binom{N_X}{2k-2\ell}}\, \frac{1}{f_0} -n^{1/2}\log
n \ge (1-2\e)\n_{\ell}p_0\,,
$$
due to our assumption on $\e$.

Now we treat typical co-degrees in the graph $G_i$. First fix $b_1,b_2\in
B_i$ and and $Y_{i,b_1},Y_{i,b_2}$ and expose a random set $X_i$. 
Let $Z_{b_1,b_2}^*$ be the number of subsets $S_1\subset
[n]$ of cardinality $|S_1|=2(k-\ell)$ such that $S_1\subset X_i$ and
both $S_1\cup Y_{i,b_1}$ and $S_1\cup Y_{i,b_2}$ form an edge in
$\cE$. By our assumption $P_e$,
$$
 \E(Z_{b_1,b_2}^*)\le
(1+\e)\bfrac{k-\ell}{\ell}^{2k-2\ell}
\binom{n-2(2\ell-k)}{2(k-\ell)}p^2\,.
$$
Using assumption $P_d$ we see that changing $X_i$ by one element changes $Z_{b_1,b_2}^*$ by at most
$$\D_{b_1,b_2}=\max_{S\in \binom{[n]}{2\ell-k+1}}|N_H(S)|\leq (1+\e)\binom{n}{2(k-\ell)-1}p.$$
Applying \eqref{ranpi} we see that \qs\ for
every $b_1,b_2\in B_i$,
\begin{equation}\label{Zb1b2*}
Z_{b_1,b_2}^*\le
\brac{1+\frac{3}{2}\e}\bfrac{k-\ell}{\ell}^{2(k-\ell)}
\binom{n}{2(k-\ell)}p^2\,.
\end{equation}

Assume $X_i$ is chosen so that (\ref{Zb1b2*}) holds. Expose the
random permutation $\s_i$ of $X_i$. Let $Z_{b_1,b_2}$ be the
co-degree of $b_1,b_2$ in $G_i$, which is the number of blocks
$X_{i,a}$ of $X_i$ under $\s_i$ such that $E_1=X_{i,a}\cup
X_{i,a+1}\cup Y_{i,b_1},E_2=X_{i,a}\cup X_{i,a+1}\cup Y_{i,b_2}\in
\cE$, and both edges $E_1,E_2$ are labeled by $i$. Then, recalling
that an edge $E\in \cE$ is labeled by $i$ with probability
$\frac{1}{f(E)}\le\frac{1}{\big(1-\frac{3}{2}\e\big)f_0}$, we get
$$
\E(Z_{b_1,b_2}) \le Z_{b_1,b_2}^*
\frac{\n_{\ell}}{\binom{N_X}{2k-2\ell}}\,
\frac{1}{\big(1-\frac{3}{2}\e\big)^2f_0^2}\ .
$$
Transposing one pair of elements of $\s_i$ changes $Z_{b_1,b_2}$ by at most 4. Using \eqref{ranpi} again, we
obtain that \qs\ for every partition $i$ and every pair $b_1,b_2\in
B_i$, the co-degree $Z_{b_1,b_2}$ of $b_1,b_2$ in $G_i$ satisfies:
\begin{eqnarray*}
Z_{b_1,b_2}&\le&
\brac{1+\frac{3}{2}\e}\bfrac{k-\ell}{\ell}^{2k-2\ell}
\binom{n}{2k-2\ell}p^2\frac{\n_{\ell}}{\binom{N_X}{2k-2\ell}}\,
\frac{1}{\big(1-\frac{3}{2}\e\big)^2f_0^2}-n^{1/2}\log n\\
&\le&(1+5\e)\n_{\ell}p_0^2\ .
\end{eqnarray*}

Now consider $a_1,a_2\in
A_i$ and and $X_{i,a_1},X_{i,a_1+1},X_{i,b_2},X_{i,b_2+1}$ and expose a random set $Y_i$. 
Let $Z_{a_1,a_2}^*$ be the number of subsets $S_1\subset
[n]$ of cardinality $|S_1|=2\ell-k$ such that $S_1\subset Y_i$ and
both $S_1\cup X_{i,a_1}\cup X_{i,a_1+1}$ and $S_1\cup X_{i,a_2}\cup X_{i,a_2+1}$ form an edge in
$\cE$. By our assumption $P_f$,
$$
 \E(Z_{a_1,a_2}^*)\le
(1+\e)\bfrac{2\ell-k}{\ell}^{2\ell-k}
\binom{n}{2\ell-k}p^2\,.
$$
Using assumption $P_c$ we see that changing $Y_i$ by one element changes $Z_{a_1,a_2}^*$ by at most
$$\D_{a_1,a_2}=\max_{S\in \binom{[n]}{2(k-\ell)+1}}|N_H(S)|\leq (1+\e)\binom{n}{2\ell-k-1}p.$$
Applying \eqref{ranpi} we see that \qs\ for
every $b_1,b_2\in B_i$,
\begin{equation}\label{Zb1b2**}
Z_{a_1,a_2}^*\le
\brac{1+\frac{3}{2}\e}\bfrac{k-\ell}{\ell}^{2k-2\ell}
\binom{n}{2k-2\ell}p^2\,.
\end{equation}

Assume $Y_i$ is chosen so that (\ref{Zb1b2**}) holds. Expose the
random permutation $\t_i$ of $Y_i$. Let $Z_{a_1,a_2}$ be the
co-degree of $a_1,a_2$ in $G_i$, which is the number of blocks
$Y_{i,b}$ of $Y_i$ under $\t_i$ such that $E_1=X_{i,a_1}\cup
X_{i,a_1+1}\cup Y_{i,b},E_2=X_{i,a_2}\cup X_{i,a_2+1}\cup Y_{i,b}\in
\cE$, and both edges $E_1,E_2$ are labeled by $i$. Then, recalling
that an edge $E\in \cE$ is labeled by $i$ with probability
$\frac{1}{f(E)}\le\frac{1}{\big(1-\frac{3}{2}\e\big)f_0}$, we get
$$
\E(Z_{a_1,a_2}) \le Z_{a_1,a_2}^*
\frac{\n_{\ell}}{\binom{N_Y}{2\ell-k}}\,
\frac{1}{\big(1-\frac{3}{2}\e\big)^2f_0^2}\ .
$$
Transposing one pair of elements of $\t_i$ changes $Z_{a_1,a_2}$ by at most 4. Using \eqref{ranpi} again, we
obtain that \qs\ for every partition $i$ and every pair $a_1,a_2\in
B_i$, the co-degree $Z_{a_1,a_2}$ of $a_1,a_2$ in $G_i$ satisfies:
\begin{eqnarray*}
Z_{a_1,a_2}&\le&
\brac{1+\frac{3}{2}\e}\bfrac{2\ell-k}{\ell}^{2\ell-k}
\binom{n}{2\ell-k}p^2\frac{\n_{\ell}}{\binom{N_Y}{2\ell-k}}\,
\frac{1}{\big(1-\frac{3}{2}\e\big)^2f_0^2}-n^{1/2}\log n\\
&\le&(1+5\e)\n_{\ell}p_0^2\ .
\end{eqnarray*}

 \proofend

\medskip

We can now apply Lemma \ref{lem2} below with
$N=\n_\ell,d=p_0,\th=5\e$ to show that each $G_{i}$ contains at
least $(1-(5\e)^{1/3})\n_\ell p_0\geq (1-2\e^{1/3})m/rn$ edge
disjoint perfect matchings. This will complete the proof of Theorem
\ref{th1} for the case $\ell<k<2\ell$.

\begin{lemma}\label{lem2}
Let $G$ be a bipartite graph with vertex set $A\cup B$ where
$|A|=|B|=N$. Suppose that the minimum degree in $G$ is at least
$(1-\th)dN$ and the maximum co-degree of two vertices is at
most $(1+\th)d^2N$ for some small value $\th\ll 1$. Suppose further
that $\th^{4/3}d^2 N\gg 1$. Then $G$ contains a collection of
$(1-\th^{1/3})dN$ edge disjoint perfect matchings.
\end{lemma}

The assumption $\th^{4/3}d^2 N\gg 1$ in the above lemma is mostly
for convenience and is implied in our circumstances  by the
assumption $\e^5\gg \log^3n/(n^{1/2}p^2)$ of Theorem \ref{th1}; it can
be relaxed somewhat.

\proofstart Let $d_0=(1-\th)d$ and $d_1=(1-\th^{1/3})d$. Going back
to \eqref{1} we see that we need to show that
\begin{equation}\label{1a}
m\geq (k-\ell)d_1N
\end{equation}
for all $K\subseteq A,L\subseteq B$,
$|K|=k,|L|=\ell,m=e(K,B\setminus L)$. Obviously we can assume
$k>\ell$.

Now,
$$m\geq k(d_0N-\ell)$$
and so \eqref{1a} is satisfied if
$$k(d_0N-\ell)\geq (k-\ell)d_1N$$
or
$$k(d_0-d_1)N\geq (k-d_1N)\ell.$$
In particular, \eqref{1a} holds if $\ell\leq (d_0-d_1)N$. Furthermore, we also have
$$
m\geq (N-\ell)(d_0N-(N-k)).
$$
If $N-\ell\leq (d_0-d_1)N$ then this implies that \eqref{1a}
holds.

So we assume from now on that
\beq{kl}
(d_0-d_1)N< \min\set{\ell,N-\ell}.
\eeq
We can further assume that $\ell\leq N/2$. For $\ell>N/2$ we can reverse the roles of $A,B$ and show
that $e(B\setminus L,A\setminus(A\setminus K))\geq ((N-\ell)-(N-k))d_1N$, which is \eqref{1a}.

We now perform the usual double counting trick by estimating the
number of paths of the form $K,B,K$ in two ways. On one hand, each
such path corresponds to a common neighbor of a pair of vertices
$a_1,a_2\in K$. Therefore, the quantity to be estimated is at most
$\binom{k}{2}d_2^2N$, where $d_2=(1+\th)^{1/2}d$. On the other hand,
it is exactly
$$
\sum_{b\in B} \binom{d(b,K)}{2}=\sum_{b\in B\setminus
L}\binom{d(b,K)}{2}+\sum_{b\in L}\binom{d(b,K)}{2}\,
$$
where $d(b,K)$ is the number of neighbors of $b$ in $K$ in the graph
$G$. Since $\sum_{b\in B\setminus L}d(b,K)=m$, we can estimate the
first summand as follows:
$$
\sum_{b\in B\setminus L}\binom{d(b,K)}{2}\ge
(N-\ell)\binom{\frac{m}{N-\ell}}{2}=\frac{m\left(\frac{m}{N-\ell}-1\right)}{2}
\ .
$$
As for the second summand, the number of edges between $L$ and $K$
can be estimated from below by $d_0Nk-m$, and therefore
$$
\sum_{b\in L}\binom{d(b,K)}{2}\ge
\ell\binom{\frac{d_0kN-m}{\ell}}{2}=
\frac{d_0kN-m}{2}\left(\frac{d_0kN-m}{\ell}-1\right)\ .
$$
It follows that
$$
m\left(\frac{m}{N-\ell}-1\right)+(d_0kN-m)\left(\frac{d_0kN-m}{\ell}-1\right)\le
k(k-1)d_2^2N\,.
$$
After performing straightforward arithmetic manipulations, we get
to:
$$
(m-d_0k(N-\ell))^2\le k(N-\ell)(k\ell d_2^2-k\ell d_0^2+d_0\ell)\,.
$$
Recalling the definitions of $d_0$ and $d_2$, we see that
$$
d_2^2-d_0^2=(1+\th-(1-\th)^2)d^2=(3\th+\th^2)d^2\,.
$$
Also, since $k\ge \ell$ and $\ell\ge (d_0-d_1)N\ge \th^{1/3}dN/2$ by
(\ref{kl}), we see that $\th dk\ge \th^{4/3}d^2N/2\gg 1$ by the
lemma's assumption. Hence $d_0\ell \ll \th d^2k\ell$. We thus arrive
at the following inequality:
\begin{equation}\label{esm-m}
(m-d_0k(N-\ell))^2\le 4\th d^2k^2\ell(N-\ell)\ .
\end{equation}
Since $\ell\le N/2$, we have 
\begin{eqnarray*}
m&\geq& (N-\ell)kd_0\brac{1-\frac{2\ell^{1/2}\th^{1/2}}{(1-\th)(N-\ell)^{1/2}}}\\
&\geq&(N-\ell)kd_0\brac{1-\frac{2\th^{1/2}}{1-\th}}.
\end{eqnarray*}
This implies \eqref{1a} if
$$\frac{d_1}{d}\leq \frac{(N-\ell)k}{(k-\ell)N}\brac{1-\frac{2\th^{1/2}}{1-\th}}(1-\th).$$
Since $\frac{k(N-\ell)}{(k-\ell)N}\ge 1$, it is enough to verify
that
$$\frac{d_1}{d}\leq
\brac{1-\th-2\th^{1/2}}.$$ 
This
is implied by
$$d_1\leq d(1-\th^{1/3}).$$
Thus there will be $d(1-\th ^{1/3})$ edge disjoint perfect
matchings. \proofend

{\bf Case 2:} $k=2\ell$ (Theorem \ref{th3}).

When $\ell=k/2$ we have $N_Y=0$ and the argument above breaks down. We can however use our result from \cite{FK1}
to obtain something.

We will construct the Hamilton cycles via the following algorithm: We first chhose $f_0$ such that
$$\frac{\log n}{\e^2}\ll f_0 \ll \min\set{\frac{\e^2np}{\log n},\frac{\e^3np}{\log n}}.$$
\begin{description}
\item[$F_1$:] Let
$$r=\frac{(2-\e)\binom{n}{k}f_0}{\n_{\ell}^2};\ \
p_0=\frac{p}{f_0}.$$
\item[$F_2$:]
Now choose $r$ random partitions $\cP_i,i=1,2,\ldots r$, of $[n]$
into $\n_\ell$ sets of size $\ell$.
\item[$F_3$:] For each $E\in \cE$ we let $f(E)$ denote the number of $i$ such that $\cP_i$ contains a
a pair of parts $X,Y$ such that $X\cup Y=E$. The random variable
$f(E)$ is distributed as $Bin(r,\r)$ where
$$\r=\frac{\binom{\n_{\ell}}{2}}{\binom{n}{k}}\,.$$
So, \qs\ $(1-\e)f_0\leq
f(E)\leq f_0$. Choose one of these $f(E)$ instances at random and
label the edge $E$ with the chosen $i$.
\item[$F_4$:] For each $i$ let $G_{i}$ be the graph obtained from $G_{\cP_i}$
by including only edges with label $i$. We will show below in Lemma
\ref{regu} that \qs\ each $G_i$ is $\brac{(1-2\e)p_0,2\e
p_0}$-regular.
\item[$F_5$:] We then apply Theorem \ref{FK1} to show that \qs\ each $G_i$ contains at least
$\brac{(1-2\e)\frac{p_0}{2}-8\e p_0}\n_\ell$ edge disjoint Hamilton
cycles. Each such Hamilton cycle corresponds to a Hamilton cycle of
type $\ell$ in $H$, and the so obtained Hamilton cycles in $H$ are
edge disjoint.
\end{description}
Thus $H$ contains at least
$$r\brac{(1-2\e)\frac{p_0}{2}-8\e p_0}\n_\ell\geq \frac{\binom{n}{k}p}{n_\ell}(1-20\e)$$
edge disjoint type $\ell$ Hamilton cycles, proving Theorem \ref{th3}.

\begin{lemma}\label{regu}
Let $G$ be a $\n$ vertex, $(\a,\e)$-regular graph. Suppose that
$G_0$ is the random subgraph of $G$ where each edge $e$ of $G$ is
included independently with probability $p_e$, where $(1-\e)p^*\leq
p_e\leq p^*=\Theta(1/\log^2\n)$. Suppose that
$$\e^2\a\n p^*\gg \log \n\ and\ \a\e^3\n p^*\gg \log 1/\e.$$
Then $G_0$ is $((1-2\e)\a p^*,2\e\a p^*)$-regular, \qs.
\end{lemma}
\proofstart The degree of vertex $v$ in $G_0$ dominates
$Bin(\a\n,(1-\e)p^*)$ and so Property $Q_a$ holds from Chernoff
bounds. Similarly the number of edges between two sets $S,T$
dominates $Bin((\a-\e)|S|\,|T|,(1-\e)p^*)$ and is dominated by
$Bin((\a-\e)|S|\,|T|,p^*)$ and Property $Q_b$ also holds from
Chernoff bounds. \proofend

{\bf Case 3:} $k=\ell$ (Theorem \ref{th4}).

Here the aim is to find many edge disjoint perfect matchings. We
construct them via the following algorithm: We first choose $f_0$ such that \eqref{f0l} holds.
\begin{description}
\item[$G_1$:] Let
$$r=(1-\e)\binom{n}{k}\frac{f_0}{\n_k^2};\ \
p_0=\frac{p}{f_0}.$$
\item[$G_2$:] Let $k_X=\rdown{k/2}$ and $k_Y=\rdup{k/2}$ and
$$N_X=\frac{k_X}{k}n\ and\ N_Y=\frac{k_Y}{k}n.$$
\item[$G_3$:]
Choose $r$ random partitions $(X_i,Y_i),i=1,2,\ldots r$, of $[n]$
into two sets of size $N_X$ and $N_Y$ respectively.

We use the notation
$$X_i=\set{x_{i,1}<x_{i,2}<\cdots<x_{i,N_X}}\ and\ Y_i=\set{y_{i,1}<y_{i,2}<\cdots<y_{i,N_Y}}.$$
\item[$G_4$:] For each $i$ we let $\s_i$ be a random permutation of $X_i$ and let
$\t_i$ be a random permutation of $Y_i$. Form the partition
$X_{i,a},a=1,2,\ldots,\n_k$, of $X_i$ into sets of size $k_X$ and
the partition $Y_{i,b},b=1,2,\ldots,\n_k$, of $Y_i$ into sets of
size $k_Y$. Here
$X_{i,a}=\set{x_{i,\s_i((a-1)k_X+1)},\ldots,x_{i,\s_i(ak_X)}}$ and
$Y_{i,b}=\set{y_{i,\t_i((b-1)k_Y+1)},\ldots,y_{i,\t_i(bk_Y)}}$.

\item[$G_5$:] Suppose now that for $E\in \cE$ there are $f(E)$ instances
$i$ such that for some $a,b$ and some partition $S_1,S_2$ of $E$ we
have $S_1=X_{i,a}$ and $S_2=Y_{i,b}$. We say that $i$ includes $E$.
Choose one of the $f(E)$ instances at random and label edge $E$ with
the chosen $i$.

Thus $f(E)$ is distributed as $Bin(r,\r)$ where
$$\r=\frac{\n_k^2}{\binom{n}{k}}.$$
So, $\big(1-\frac{3}{2}\e\big)f_0\leq f(E)\leq f_0$ \qs.
\item[$G_6$:] Let $G_{i}$ be the bipartite graph with vertex partition $A_{i}$ and
$B_i$ comprising disjoint copies of $[\n_k]$. For $a\in A_{i}$ and
$b\in B_i$  we make $(a,b)$ an edge of $G_{i}$ if $E=X_{i,a}\cup
Y_{i,b}\in \cE$ and $E$ is labelled with $i$. So, by construction,
each $e\in E$ is associated with at most one $G_{i}$.
\item[$G_7$:] We claim (see Lemmas \ref{lem2} and \ref{lem1z}) that \whp\
each $G_{i}$ will contain at least
$$n_0=(1-2\e^{1/3})\n_k p_0$$
edge disjoint
perfect matchings.
\end{description}
So $H$ will contain at least $rn_0$ edge disjoint perfect matchings,
completing the proof of Theorem~\ref{th4}.
\begin{lemma}\label{lem1z}
{\bf Whp}, over our random choices of $X_i,Y_i,\s_i,\t_i$, each
$G_{i}$ has minimum degree at least $(1-2\e)\n_k p_0$ and maximum
co-degree at most $(1+5\e)\n_k p_0^2$.
\end{lemma}
\proofstart The arguments here are very similar to those in Lemma
\ref{lem3}, so we will be rather brief. We fix $i$ and focus on
$G_{i}$. We first show that the minimum degree in $G_{i}$ is large.
For $a\in A_{i}$, denote by $Z_a$ its degree in $G_i$. Then, using
assumptions $R_a$ and $R_c$ and martingale-type arguments, we can
show that
$$
\E(Z_a)\ge
\big(1-\frac{3}{2}\e\big)\bfrac{k_Y}{k}^{k_Y}\binom{n}{k_Y}p\,
\frac{\n_{k}}{\binom{N_Y}{k_Y}}\, \frac{1}{f_0}\ .
$$
Using concentration results for permutation graphs again, we derive
that \qs\ for every partition $i$ and every $a\in A_i$, the degree
of $a$ in $G_i$ is at least
$$
\big(1-\frac{3}{2}\e\big)\bfrac{k_Y}{k}^{k_Y}\binom{n}{k_Y}p\,
\frac{\n_{k}}{\binom{N_Y}{k_Y}}\, \frac{1}{f_0}-n^{1/2}\log n\ge
(1-2\e)\n_kp_0\,,
$$
due to our assumption on $\e$.

Let now $Z_b$ denote the degree of vertex $b\in B_i$ in $G_i$. We
can argue similarly, while invoking assumptions $R_b$, $R_d$, to
show that \qs\ for every partition $i$ and every $b\in B_i$,
$$
Z_b\ge (1-2\e)\n_kp_0\,.
$$

Finally, we treat the co-degrees of pairs of vertices in $G_i$.
Let $b_1,b_2\in B_i$.
Let $Z_{b_1,b_2}$ be their co-degree in $G_i$. Then using assumption
$R_e$ and martingale-type concentration arguments, we can prove that
\qs\ for every partition $i$ and every pair of vertices $b_1,b_2\in
B_i$
\begin{eqnarray*}
Z_{b_1,b_2}&\le& \big(1+\frac{3}{2}\e\big)\bfrac{k_X}{k}^{k_X}
\binom{n}{k_X}p^2\frac{\n_{k}}{\binom{N_X}{k_X}}\,
\frac{1}{\big(1-\frac{3}{2}\e\big)^2f_0^2}-n^{1/2}\log n\\
&\le&(1+5\e)\n_{k}p_0^2\ .
\end{eqnarray*}
Similarly, if $a_1,a_2\in A_i$,
let $Z_{a_1,a_2}$ be their co-degree in $G_i$. Then using assumption
$R_f$ and martingale-type concentration arguments, we can prove that
\qs\ for every partition $i$ and every pair of vertices $a_1,a_2\in
A_i$
\begin{eqnarray*}
Z_{a_1,a_2}&\le& \big(1+\frac{3}{2}\e\big)\bfrac{k_Y}{k}^{k_Y}
\binom{n}{k_Y}p^2\frac{\n_{k}}{\binom{N_Y}{k_Y}}\,
\frac{1}{\big(1-\frac{3}{2}\e\big)^2f_0^2}-n^{1/2}\log n\\
&\le&(1+5\e)\n_{k}p_0^2\ .
\end{eqnarray*}
\proofend

\section{Concluding remarks}
In this paper we have derived sufficient conditions for packing
almost edges of $k$-uniform random and pseudo-random hypergraphs
into disjoint type $\ell$ Hamilton cycles. This appears to be a
first result of this kind for the problem of packing Hamilton cycles
in this setting. There is no reason to believe our assumptions on
the edge probability $p(n)$ or the density of a pseudo-random
hypergraph are tight, and it would be quite natural to try and
extend them and to obtain tight(er) bounds.

Our results say nothing at all for the case $\ell <k/2$. It would be
nice to extend our results to this, apparently harder, case.

In our paper we managed to circumvent the absence of
results and techniques for the appearance of a Hamilton cycle of
essentially any type in a random hypergraph $H(n,p,k)$ by going to
much larger edge probabilities/densities than any plausible guess
for the threshold for the appearance of a single Hamilton cycle. It
would be very interesting to address specifically the question of
the appearance of a Hamilton cycle in the random hypergraph and to
derive better upper bounds on the corresponding threshold.

In the paper \cite{Jumble} we were able to show how to use the results of 
\cite{FK1} in a game theoretic setting. More precisely, we showed how to 
play a Maker-Breaker type of game on the complete graph where Maker 
is able to construct an $(1/2-\e,\e)$-regular graph, $\e=o(1)$. We could then
use the results of \cite{FK1} to show that Maker could construct approximately $n/4$
edge disjoint Hamilton cycles when alternately choosing edges against an adversary. 
The techniques of that paper can be extended to the hypergraph setting in a straightforward manner.


\begin{thebibliography}{99}
\bibitem{BS} R. F. Bailey and B. Stevens, {\em Hamiltonian
decompostions of complete $k$-uniform hypergraphs}, Discrete
Mathematics, in press.
\bibitem{Bar}
Zs. Baranyai, On the factorization of the complete uniform
hypergraph, In: {\em Infinite and finite sets}, Vol. I, Colloq.
Math. Soc. Janos Bolyai, Vol. 10, North-Holland, Amsterdam, 1975,
91--108.
\bibitem{Bol}
 B. Bollob\'as, {\bf Random graphs}, 2nd ed., Vol. 73, Cambridge
University Press, Cambridge, 2001.
\bibitem{CGW}
F. R. K. Chung, R. L. Graham, and R. M. Wilson, {\em Quasi-random
graphs}, Combinatorica 9 (1989), 345--362.
\bibitem{FK1} A.M. Frieze and M. Krivelevich, {\em On packing Hamilton cycles in $\e$-regular graphs},
Journal of Combinatorial Theory Ser. B 94 (2005) 159-172.
\bibitem{Jumble}
A.M. Frieze, M. Krivelevich, O. Pikhurko and T. Szabo, {\em The game
of JumbleG}, Combinatorics, Probability and Computing 14 (2005)
783-794.
\bibitem{FP} A.M. Frieze and B. Pittel, Perfect matchings in random graphs with prescribed minimal degree,
{\em Trends in Mathematics, Birkhauser Verlag, Basel} (2004) 95-132.
\bibitem{HS} H. H\'an and M. Schacht, {\em Dirac-type results for
loose Hamilton cycles in uniform hypergraphs}, Journal of
Combinatorial Theory  Ser. B 100 (2010), 332-346.
\bibitem{JLR}
 S. Janson, T. \L uczak and A. Ruci\'nski, {\bf Random
graphs},  Wiley, New York, 2000.
\bibitem{JKV}
A. Johansson, J. Kahn and V. Vu, {\em Factors in random graphs},
Random Structures and Algorithms 33 (2008), 1--28.
\bibitem{KKMO}
P. Keevash, D. K\"uhn, R. Mycroft and D. Osthus, {\em Loose Hamilton
cycles in hypergraphs}, manuscript.
\bibitem{KS}
M. Krivelevich and B. Sudakov, Pseudo-random graphs, In: {\em More
sets, graphs and numbers}, E. Gy\"ori, G. O. H. Katona, L. Lov\'asz,
Eds., Bolyai Soc. Math. Studies Vol. 15, 199--262.
\bibitem{KMO}
D. K\"uhn, R. Mycroft and D. Osthus, {\em Hamilton $\ell$-cycles in
uniform hypergraphs}, manuscript.
\bibitem{McD} C. McDiarmid, {Concentration}, in:
{\em Probabilistic Methods for Algorithmic Discrete Mathematics},
Algorithms Combin., 16, Springer, Berlin, 1998, 195--248.
\bibitem{Thom}
A. Thomason, Pseudorandom graphs. In: {\em Random graphs '85
(Pozna\'n, 1985)}, North-Holland Math. Stud., 144, North-Holland,
Amsterdam, 1987, 307--331.

\end{thebibliography}
\end{document}